\documentclass{amsart}
\usepackage[latin1]{inputenc}
\usepackage{amssymb}
\usepackage{amsmath}
\usepackage{epsfig}

\newtheorem{theorem}{Theorem}[section]
\newtheorem{rem}{Remark}[section]
\newtheorem{ex}{Example}[section]
\newtheorem{prop}{Proposition}[section]
\newtheorem{lemma}{Lemma}[section]

\newtheorem{defi}{Definition}[section]

\newtheorem{quest}{Question}[section]

\newcommand{\Dh}{\mathcal{D}}
\newcommand{\lra}{\longrightarrow}

\newcommand{\bprf}{{\it Proof.~}}

\newcommand{\ra}{\rightarrow}
\newcommand{\eprf}{\hfill $\square$ \smallskip\par}

\newcommand{\PP}{ \mathbb{P}}
\newcommand{\C }{ \mathbb{C}}

\newcommand{\Z}{\mathbb{Z}}
\newcommand{\Q}{\mathbb{Q}}

\DeclareMathOperator{\rank}{\rm{rank}}

\DeclareMathOperator{\Fix}{\rm{Fix}}

\makeatletter
\def\blfootnote{\xdef\@thefnmark{}\@footnotetext}

\begin{document}

\title{On symplectic and non--symplectic automorphisms of K3 surfaces}

\author{Alice Garbagnati and Alessandra Sarti}
\address{Alice Garbagnati, Dipartimento di Matematica, Universit\`a di Milano,
  via Saldini 50, I-20133 Milano, Italy}

\email{alice.garbagnati@unimi.it}
\urladdr{http://sites.google.com/site/alicegarbagnati/home}

\address{Alessandra Sarti, Universit\'e de Poitiers, 
Laboratoire de Math\'ematiques et Applications, 
 T\'el\'eport 2 
Boulevard Marie et Pierre Curie
 BP 30179,
86962 Futuroscope Chasseneuil Cedex, France}

\email{sarti@math.univ-poitiers.fr}
\urladdr{http://www-math.sp2mi.univ-poitiers.fr/~sarti/}

\begin{abstract}
In this paper we investigate when the generic member of a family of K3 surfaces admitting a non--symplectic automorphism of finite order admits also a symplectic automorphism of the same order. We give a complete answer to this question if the order of the automorphism is a prime number and we provide several examples and partial results otherwise. Moreover we prove that, under certain conditions, a K3 surface admitting a non--symplectic automorphism of prime odd order, $p$, also admits a non--symplectic automorphism of order $2p$. This generalizes a previous result by J. Dillies for $p=3$.
\end{abstract}

\subjclass[2000]{Primary 14J28; Secondary 14J50, 14J10}
\keywords{$K3$ surfaces, automorphisms}

\maketitle

\pagestyle{myheadings}
\markboth{Alice Garbagnati and Alessandra Sarti}{On symplectic and non--symplectic automorphisms of K3 surfaces}

\section{Introduction}
 An automorphism of finite order $n$ on a K3 surface is called {\it symplectic} if it acts trivially on the holomorphic 2--form of the K3 surface and it is called {\it purely non--symplectic} if it acts as a multiplication by a primitive $n$--th root of the unity. These notions were introduced by V.V. Nikulin in the 80's (cf. \cite{niksympl}). \\
In \cite{niksympl} the finite abelian groups $G$  which act symplectically on a K3 surface are classified and it is proved that the existence of a primitive embedding of a certain negative definite lattice $\Omega_G$ (depending only on the group $G$) in the N\'eron--Severi group of a K3 surface $X$ is equivalent to the fact that $G$ acts symplectically on $X$. In \cite{aliale1}, \cite{aliale2} the lattice $\Omega_G$ is computed for each finite abelian group $G$ and, thanks to the Nikulin's result, the families of K3 surfaces admitting $G$ as group of symplectic automorphisms are described as families of $L_G$-polarized K3 surfaces, for certain lattices $L_G$.\\
On the other hand in \cite{niknonsympl}, \cite{michiale}, \cite{michialet} families of K3 surfaces admitting a non--symplectic automorphism of prime order and a certain fixed locus  are classified. These are also $M$-polarized K3 surfaces, for a certain lattice $M$ depending on the order and on the fixed locus of the automorphism.\\
One of the aim of this paper is to bring these two descriptions together. \\

In \cite{alidih} and \cite{Dillies} a particular phenomenon is described for certain groups of symplectic and non--symplectic automorphisms respectively: each K3 surface admitting a particular group $G$ as group of symplectic (resp. non--symplectic) automorphisms, automatically admits a larger group $H$ of symplectic (resp. non--symplectic) automorphisms (see Theorems \ref{theorem: Z5 and D5},\ref{theorem: non--symplectic 3 iff 6}). In particular this implies that the family of the K3 surfaces admitting $G$ as group of (non--)symplectic automorphisms coincides with the family of the K3 surfaces admitting $H$ as group of (non--)symplectic automorphisms.\\ 
Due to these results it seems natural to ask:
\begin{quest}\label{question} When does the generic member of a family of K3 surfaces admitting a certain group $G$ of automorphisms admit a larger group $H$ of automorphisms? \end{quest}
In Section \ref{section: automorphisms on K3 surfaces} we generalize the results of Dillies (cf. \cite{Dillies}) assuming that the groups $G$ and $H$ both act non--symplectically on the K3 surface. More precisely, let $p$ be an odd prime number. Then, under some conditions on the fixed locus of the non--symplectic automorphism of order $p$, we prove that a K3 surface admits $\Z/p\Z$ as group of non--symplectic automorphisms if and only if it admits  $\Z/2p\Z$ as group of non--symplectic automorphisms (cf. Theroem \ref{theorem: non symplectic p iff 2p}). There are two exceptions, in fact if the fixed locus of the automorphism of order $p$ consists of isolated points and $p=7,11$, then the same result is false (cf. Theorem \ref{711}).\\
We analyze newt the case when $G=<\eta>$ is a finite cyclic non--symplectic group and $H$ is generated by $\eta$ and by a symplectic automorphism $\sigma$. Very much is known about the groups $\langle \eta\rangle$ and $\langle \sigma\rangle$ and it turns out that the order  $o(\sigma)\leq 8$ and $o(\eta)\leq 66$. Thus the order $2\leq n\leq 8$ is possible for both $\eta$ and $\sigma$. One can hence ask if there are K3 surfaces having both the automorphisms $\sigma$ and $\eta$ with $o(\sigma)=o(\eta)$. The answer is positive and surprising, in the sense that in some cases the {\it generic} K3 surface with  non--symplectic automorphism also admits a symplectic automorphism of the same order. The main result of this paper is summarized in the following: 
\begin{theorem}\label{theorem: main results} Assume that $X$ is a generic K3 surface with a non--symplectic automorphism $\eta$ of order $p$ then :
\begin{itemize}
\item If $p=2$ then $X=X_{r,a,\delta}$ (the Nikulin invariants $r,a,\delta$ will be introduced in Section \ref{section: K3 surfaecs with sympl and non sympl}) admits also a symplectic involution if and only if either $\delta=1$ and $a>16-r$ or $\delta=0$ and $a>16-r$ or $\delta=0$ and $a=6,\, r=10$. 
\item If $p=3$, then $X$ admits a symplectic automorphism of order 3 if and only if the fixed locus of $\eta$ consists of $n$ points and $n-3$ curves with $n\geq 6$.
\item If $p=5,7$ then $X$ does not admit a symplectic automorphism of the same order (in particular, if $p=7$ such an automorphism does not exist for any K3 surface with a non symplectic automorphism of order 7, not only for the generic one in a family).
\end{itemize}
\end{theorem}
We give the proof in the Theorems \ref{theorem: no sympl e non for p=7,8}, \ref{theorem: syml and non iff, order 2}, \ref{theorem: order 3},  \ref{theorem: order 5 non sympl and non}. We remark that the theorem does not say anything about existence of special K3 surfaces with a symplectic and a non-symplectic automorphism of the same order (at least if $p\neq 7$), indeed in Example \ref{ex: rigid order 5} we describe a rigid K3 surface admitting both a symplectic and a non symplectic automorphism of order 5.\\
We obtain complete results only in the cases that the order of the automorphism is a prime number. Indeed the classification of families of K3 surfaces with non--symplectic automorphisms of non prime order is not complete. However we show that there exist no K3 surfaces admitting both a symplectic and a non--symplectic automorphism of order 8 (cf. Table \ref{tablerank}) and we show that there exists a 1-dimensional family of K3 surfaces admitting both a symplectic and a non--symplectic automorphisms of order 4 (resp.  6) (cf. Example \ref{ex: 1 dim order 4}, resp. Example \ref{ex: 1 dim order 6}). Moreover, in case of order 6, we show that the generic member of a family of K3 surfaces admitting a non--symplectic automorphism of order 6, does not admit a symplectic automorphism of the same order (cf.\ Theorem \ref{theorem: order 3 and 6 sympl}). The same result is proved for the order 4 if the dimension of the family is at least 2 (cf. Theorem \ref{theorem: no 2 dimensional family order 4}).\\

{\bf Acknowledgements.} We warmly thank Bert van Geemen and Michela Artebani for useful comments and discussions. During the preparation of the paper the first author was visiting several times the University of Poitiers  and the second author was visiting several times the University of Milan, they benefitted partially of a financial support by {\it GRIFGA : Groupement de Recherche europ\'een Italo-Fran\c cais en G\'eom\'etrie Algébrique}. The paper was completed during the second author stay at the Institut Henri Poincar\'e in Paris for the {\it Trimestre de G\'eom\'etrie Alg\'ebrique}, she would like to thank this institution for support and for the stimulating working atmosphere.\\

\section{Automorphisms on K3 surfaces}\label{section: automorphisms on K3 surfaces}

\begin{defi} Let $X$ be a smooth complex surface. The surface $X$ is a K3 surface if the canonical bundle of $X$ is trivial and the irregularity of $X$, $q(X):=h^{1,0}(X)$, is zero.\\
If $X$ is a K3 surface, then $h^{2,0}(X)=1$. We choose a generator $\omega_X$ of $H^{2,0}(X)$, which is called the period of $X$.
\end{defi}
The second cohomology group of a K3 surface, equipped with the cup product, is isometric to a lattice, which is the unique, up to isometries, even unimodular lattice with signature $(3,19)$. This lattice will be denoted by $\Lambda_{K3}$ and is isometric to $U\oplus U\oplus U\oplus E_8\oplus E_8$, where $U$ is the unimodular lattice with bilinear form $\left[\begin{array}{ll}0&1\\1&0\end{array}\right]$ and $E_8$ is the even negative definite lattice  associated to the Dynkin diagram $E_8$. The N\'eron--Severi group of a K3 surface $X$, $NS(X)$, coincides with its Picard group. The transcendental lattice of $X$, $T_X$, is the orthogonal lattice to $NS(X)$ in $H^2(X,\Z)$.  
\begin{defi} An isometry $ \alpha$ of $H^2(X,\Z)$ is an effective isometry if it preserves the K\"ahler cone of $X$. An isometry $ \alpha$ of $H^2(X,\Z)$ is a Hodge isometry if its $\C$-linear extension to $H^2(X,\C)$ preserves the Hodge decomposition of $H^2(X,\C)$. \end{defi}
\begin{theorem}{\rm (\cite{BurnsRapoporttorellitheoremK3})} Let $X$ be a K3 surface and $g$ be an automorphism of $X$, then $g^*$ is an effective Hodge isometry of $H^2(X,\Z)$. Viceversa, let $f$ be an effective Hodge isometry of $H^2(X,\Z)$, then $f$ is induced by a unique automorphism of $X$.\end{theorem}

If $g$ is an automorphism on a K3 surface $X$, then $g^*$ preserves the space of the holomorphic 2-forms on $X$, and hence $g^* (\omega_X)=\lambda\omega_X$, $\lambda\in\C^*$.

\begin{defi} An automorphism $\sigma$ on a K3 surface $X$ is symplectic if $\sigma^*$ acts as the identity on $H^{2,0}(X)$, that is $\sigma^*(\omega_X)=\omega_X$. Equivalently $\sigma$ is symplectic if the isometry induced by $\sigma^*$ on the transcendental lattice is the identity. An automorphism $ \eta$ of finite order $m$ is purely non--symplectic if $ \eta(\omega_X)=\zeta_m\omega_X$, where $\zeta_m$ is a primitive $m$-th root of unity.\end{defi}

In the following we will only say that an automorphism is non--symplectic but we mean that it is purely non--symplectic.

In \cite{niksympl} it  is proved that for each finite group $G$ of automorphisms on a K3 surface $X$, there exists the following exact sequence:
\begin{equation}\label{equ: sequence}
1\lra G_0\lra G\overset{\alpha}\lra \Gamma_m\lra 1
\end{equation}
where $\alpha$ is the natural representation of $G$ in $H^{2,0}(X)=\C\omega_X$, and $m$ is a positive
integer. Then $\Gamma_m$ is cyclic of order $m$.
If $\Gamma_m\neq\{1\}$, the K3 surface $X$ is algebraic, \cite[Theorem 3.1]{niksympl}. 
The group $\Gamma_m$ has order $m\leq 66$ and if $m$ is a prime number, then $m\leq 19$ (cf. \cite{niksympl}). A complete list of the finite groups $G_0$ acting symplectically on a K3 surface is given in \cite{Xiao} and consists of 79 groups. 
We observe that there are a priori many possibilities for the groups $G_0$ and $\Gamma_m$. However not all these possibilities correspond to different families of K3 surfaces. Indeed the following two results show that requiring a certain group of (resp.\ non--) symplectic automorphisms on a K3 surface $X$ implies that there is a bigger group of (resp.\ non--) symplectic automorphisms.
\begin{theorem}\label{theorem: Z5 and D5}{\rm (\cite{alidih})} A K3 surface admits $\Z/5\Z$ as group of symplectic automorphisms if and only if it admits the dihedral group $\mathcal{D}_5$ of order 10 as group of symplectic automorphisms.\end{theorem}
\begin{theorem}\label{theorem: non--symplectic 3 iff 6}{\rm (\cite{Dillies})} A K3 surface admits $\Z/3\Z$ as group of   non--symplectic automorphisms if and only if it admits the  group $\Z/6\Z$ as group of   non--symplectic automorphisms.\end{theorem}
In \cite{alidih} a criterion is given which implies that a group $G$ acts symplectically on a K3 surface if and only if a larger group $H(\supset G)$ acts symplectically on it. Here we prove that also in the non--symplectic case one can extend the result of Theorem \ref{theorem: non--symplectic 3 iff 6}. 

\begin{theorem}\label{theorem: non symplectic p iff 2p} Let $p=5,13,17,19$. A K3 surface admits $\Z/p\Z$ as group of   non--symplectic automorphisms if and only if it admits $\Z/2p\Z$ as group of   non--symplectic automorphisms.\\
Let $q=7,11$. If a K3 surface $X$ admits $\Z/q\Z$ as group of   non--symplectic automorphisms and the fixed locus of this automorphism contains at least one curve, then  $X$ admits $\Z/2q\Z$ as group of   non--symplectic automorphisms.
\end{theorem}
\bprf In \cite{michialet} the general member $X$ of a family of K3 surfaces admitting a non--symplectic automorphism $ \eta$ of prime order $5\leq p\leq 19$  is described for any possible fixed locus of $\eta$. To prove the theorem, it suffices to show that this general member admits a non--symplectic automorphism of order $2p$. Let us now consider the following 2 constructions:\begin{itemize}\item[case A]
Let $\alpha_p$ be an automorphism of $\mathbb{P}^2$ of prime order $p\not=2$. Let $C_{\alpha_p}=V(f_{\alpha_p}(x_0:x_1:x_2))$ be a family of sextic plane curves invariant for $\alpha_p$. Let $S_{\alpha_p}$ be the K3 surface obtained as double cover of $\mathbb{P}^2$ branched along $C_{\alpha_p}$. A (possibly singular) model of $S_{\alpha_p}\subset W\mathbb{P}(3,1,1,1)$ is $u^2=f_{\alpha_p}(x_0:x_1:x_2)$. The automorphism acting on $W\mathbb{P}(3,1,1,1)$ as $ (u;(x_0:x_1:x_2))\mapsto (u;\alpha_p(x_0:x_1:x_2))$ restricts to an automorphism $\eta$ of $S_{\alpha_p}$. It  has order $p$. It is now clear that the surface $S_{\alpha_p}$ admits also the automorphism $\beta:(u;(x_0:x_1:x_2))\ra (-u;\alpha_p(x_0:x_1:x_2))$, which is the composition of $ \eta$ with the covering automorphism and which has order $2p$.
\item[case B]
Let $R$ be a K3 surface admitting an elliptic fibration with Weierstrass equation $y^2=x^3+A(t)x+B(t)$. If $R$ admits an automorphism of prime order $p\not=2$, $ \eta:(x,y,t)\ra (\zeta_p^a x, \zeta_p^b y,\zeta_p^c t)$, then it admits also an automorphism $\beta:(x,y,t)\ra (\zeta_p^a x, -\zeta_p^b y,\zeta_p^c t)$ of order $2p$.\end{itemize}
The general member of the families of K3 surfaces with a non--symplectic automorphism of order 5 are described in \cite{michialet} as double covers of $\mathbb{P}^2$ and the non--symplectic automorphism is induced by an automorphism of $\mathbb{P}^2$. So the non--symplectic automorphism of order 5 is constructed as $ \eta$ in case A. This shows that these K3 surfaces admit also an automorphism of order 10, costructed as $\beta$ in case A. Since the cover involution is a non--symplectic automorphism of order 2, the automorphism $\beta$ is a   non--symplectic automorphism of order 10.\\
In cases $p=13,17,19$, the general member of the families and the non--symplectic automorphism of prime order $p$ on it can be constructed as the automorphism $\eta$ in case B (cf. \cite{michialet}). The automorphism $\beta$ in case B is an automorphism of order $2p$ on these K3 surfaces. It is   non--symplectic since it is the composition of two  non--symplectic automorphisms with coprime order.
The situation is similar in cases $p=7,11$ if the automorphism fixes at least one curve. 
\eprf
\subsection{Holomorphic Lefschetz formula for a non--symplectic automorphism of order $14$ and of order $22$} The  Theorem \ref{theorem: non symplectic p iff 2p} does not describe any relation between automorphisms of order $p$ and of order $2p$ if $p=7,11$ and the non--symplectic automorphism of order $p$ fixes only isolated points. In order to describe this situation we need a deeper analysis of the non--symplectic automorphisms of order $2p$, $p=7,11$.\\
We can assume that $\eta$ acts on $\omega_X$ as the multiplication by $\zeta_m$, $m=14$ or $m=22$.
The action of $\eta$ can be locally linearized and diagonalized at a fixed point $x\in X^{\sigma}$ (see \cite{Cartan} and \cite[\S 5]{niksympl} ), so that its possible local actions are
$$ A_{m,t}=\left(\begin{array}{cc}
	\zeta_m^{t+1} &0\\
0&\zeta_m^{m-t}
\end{array}
\right),\ \ t=0,\dots,\frac{m-2}{2}.$$
If $t=0$ then $x$ belongs to a smooth fixed curve for $\eta$, otherwise $x$ is an isolated fixed point.
We will say that an isolated point $x\in X^{\eta}$ is of \emph{type }$t$ ($t>0$) if the local action at $x$ is given by $A_{m,t}$ and we will denote by $n_{t}$ the number of isolated points of $\eta$ of type $t$.

The holomorphic Lefschetz formula \cite[Theorem 4.6]{AS} allows to compute the holomorphic Lefschetz number $L(\eta )$ of $\eta$ in two ways. First we have that
\[ L(\eta ) = \sum _{i=0}^{2} (-1)^i\text{tr}(\eta ^{\ast }|H^{i}(X, \mathcal{O}_{X})). \]
Since we have $H^{2}(X, \mathcal{O}_{X})=H^{0,2}=\overline{H^{2,0}}=\overline{\mathbb{C}\omega_X}$ we obtain: 
\begin{equation}\label{l1}
L(\eta )=1+\zeta_m ^{m-1}.
\end{equation}
On the other hand, we also have that
\[ L(\eta ) = \sum_{t=1}^{m-2} n_{t} a(t)+\sum_i b(C_i), \]
where  the $C_i$ are the $\eta$-fixed curves of genus $g(C_i)$, 
\begin{equation}\label{l2}
 a(t):=\frac{1}{\det(I-\sigma^{\ast }|T_{t} )}
  =\frac{1}{\det(I-A_{m,t})} 
  =\frac{1}{(1-\zeta ^{t})(1-\zeta^{m-t+1})}, 
 \end{equation}
 with $T_{t}$ the tangent space of $X$ at a point of type $t$, and\\
\begin{equation}
b(C_i):=\frac{(1+\zeta_m)(1-g(C_i))}{(1-\zeta_m )^{2}}.
\end{equation}
Denoting by $h=\sum (1-g(C_i))$ we can then write
$$ L(\eta ) = \sum_{t=1}^{m-2} n_{t} a(t)+h\frac{(1+\zeta_m)}{(1-\zeta_m )^{2}}. $$

If $X^{\eta}$ is either empty or the union of two elliptic curves, then $L(\eta)=0$, this is possible only in case of involutions. By applying the holomorphic Lefschetz formula one obtains the following results in the case that the order of the automorphism is $14$ or $22$ and the automorphism $\eta^2$ of order $7$ or $11$ has only isolated fixed points. \\
$\fbox{m=14}$. In this case by \cite[Theorem 6.3]{michialet} the local actions at the fixed points of $\eta^2$ are
$$ \left(\begin{array}{cc}
	\zeta_7^{2} &0\\
0&\zeta_7^{6}
\end{array}
\right), \left(\begin{array}{cc}
	\zeta_7^{3} &0\\
0&\zeta_7^{5}
\end{array}
\right)\ \ 
$$

and there are two fixed points with local action of the first type and one fixed point with local action of the second type. Observe that the automorphism $\eta$ of order $14$ has only isolated fixed points too and in the first case the local action of $\eta$ can be of type $1$ or $5$ and in the second case of type $2$ or $4$. Then $\eta$ either fixes the three points fixed by $\eta^2$ or interchanges two points and fixes the third one. In any case we have $h=n_3=0$. By using the fact that the roots $\zeta_{14}^i$, $i=0,\ldots,5$ are linearly independent over $\Q$ one get a system of equations (with MAPLE): 
\begin{eqnarray*}
\left\{\begin{array}{rrrrrrrrr}
3 n_1&-& 3 n_2 &-& 5 n_4 &+& n_5&=&7\\
8 n_1 &-& n_2 &+& 3 n_4&-& 2 n_5 &=&7\\
6 n_1 &+&  n_2 &-& 3 n_4 &+& 2 n_5&=&7 \\
 4 n_1 &+ &3 n_2 &-& 9 n_4&+& 6 n_5& =&7 \\
9 n_1 &+& 5 n_2 &-& n_4 &+& 3 n_5&=&7\\
6 n_1&+& n_2 &-& 3 n_4 &+& 2 n_5 &=&7 
\end{array}\right.
\end{eqnarray*}

Subtracting the first equation from the third we get $n_5=-3n_1-4n_2-2n_4$. Since $n_i\geq 0$, $i=1,2,4,5$, we obtain  $n_1=n_2=n_4=n_5=0$, which is impossible in any other equation. Hence there are no non-symplectic automorhpisms of order 14 on K3 surfaces such that the square fixes isolated points.

$\fbox{m=22}$. In this case by \cite[Theorem 7.3]{michialet} the local actions at the fixed points of $\eta^2$  are
$$ \left(\begin{array}{cc}
	\zeta_{11}^{2} &0\\
0&\zeta_{11}^{10}
\end{array}
\right), \left(\begin{array}{cc}
	\zeta_{11}^{3} &0\\
0&\zeta_{11}^{9}
\end{array}
\right).
$$
Observe that the automorphism of order $22$ must fix also these two points.
Near to the first fixed point the local action can be of type $1$ or $9$, near to the second one is of type $4$ or $6$. In any case we have $h=n_2=n_3=n_5=n_7=n_8=n_{10}=0$.  By using the fact that the roots $\zeta_{22}^i$, $i=0,\ldots,9$  are independent over $\Q$ one get the system of equations (with MAPLE):
\begin{eqnarray*}
\left\{\begin{array}{rrrrrrrrr}
5 n_1  &-&8 n_4&+& 2 n_6& +&  n_9 &=&11  \\
24 n_1 &- & n_4&+& 3 n_6 &-&4 n_9 &=&11 \\
10 n_1 &-& 5 n_4 &+& 15 n_6 &+&2 n_9  &=&11\\
18 n_1&- &9 n_4 &+& 5 n_6&+& 8 n_9 &=& 11 \\
 15 n_1&-& 2 n_4&+& 6 n_6&+&3 n_9&=& 11\\
12 n_1&  +& 5 n_4&+&7 n_6  &-& 2 n_9&=& 11\\
20 n_1 &+ & n_4&-& 3 n_6 &+& 4 n_9  &=&11\\
6 n_1&- &3 n_4  &+& 9 n_6 &+& 10 n_9 &=& 11\\
 25 n_1& + &4 n_4  &+& 10 n_6&+& 5 n_9&=&11 \\
15 n_1& -& 2 n_4 &+& 6 n_6 &+& 3 n_9 &=& 11\\
\end{array}\right.
\end{eqnarray*}
Subtracting the 9th equation from the 7th, one get $13n_6+n_9+5n_1+3n_4=0$, which gives $n_1=n_9=n_4=n_6=0$ .
This is impossible in any other equation. Hence there are no non-symplectic automorhpisms of order 22 on K3 surfaces such that the square fixes isolated points.\\ 

Assume now that $X$ has a non-symplectic automorphism $\eta$ of order $7$ or $11$ then with the same notation as in the exact sequence \eqref{equ: sequence} $m$ is a multiple of $7$ or $11$. We show then
\begin{theorem}\label{711}
Assume that $\eta$ has only isolated fixed points.\\ 
If $X$ is a {\rm generic} K3 surface with a non--symplectic  automorphism of order seven then $G_0=\{1\}$ and $m=7$.\\
If $X$ is {\rm any} K3 surface with a non--symplectic automorphism of order eleven then $G_0=\{1\}$ and $m=11$. 
\end{theorem}
\bprf
Assume first that $X$ is generic. Then  we have that $T_X=U(7)\oplus U\oplus E_8\oplus A_6$, $NS(X)=U(7)\oplus K_7$ in the case of the order $7$  (where $K_7\simeq\left[\begin{array}{rr}-4&1\\1&-2\end{array}\right]$) and $T_X=U\oplus U(11)\oplus E_8\oplus E_8$, $NS(X)=U(11)$ in the case of the order $11$. In both cases $\rank(NS(X))<8$ so by \cite{niksympl} or \cite[Lemma 1.2]{MO} we have $G_0=\{1\}$. Since $X$ has an automorphism of order $p=7,11$, we can write $m=pn$ and one sees immediately that $n\leq 6$ (cf. \cite{niksympl}, \cite{zhang}). Let $h$ denote a generator of $G$ (since $G_0=\{1\}$, $G=\Gamma_m$ is cyclic) then $\eta=h^{n'}$ and the fixed loci $h(X^\eta)\subset X^\eta$ and $X^h\subset X^\eta$, in particular the fixed locus of $X^h$ consists of isolated points (in fact $X^h$ is not empty since $h$ is not an involution). \\
$\fbox{p=7}$. The only possible values are $m=7,14$ because only in these cases the Euler function $\varphi(m)$ divides the rank of the transcendental lattice which is $18$. We must exclude the case $m=14$. Let $X^\eta=\{P_1,P_2,P_3\}$ be the fixed points of $\eta$ on $X$. Then $X^h$ consists either of one point or of three points. By using the Lefschetz fixed point formula for an automorphism of order $14$ as done at the beginning of the section, one sees that this is not possible.\\
$\fbox{p=11}$. The possible values are $n=11,22,33,44,66$. If $n=33,44,66$, the K3 surface is unique (see \cite{kondo1}), and our family is 1-dimensional, so this is not possible. We are left to exclude the case $m=22$. Let $X^\eta=\{P_1,P_2\}$, since $X^h$ is not empty it must be equal to $X^\eta$. Again using the Lefschetz fixed point formula, here for an automorphism of order $22$ as done at the beginning of the section, one sees that this is not possible.\\
Consider now any K3 surface with non--symplectic action by an automorphism of order 11. Since the Euler function of $11$ is $10$, the rank of the transcendental lattice can be only $10$ or $20$. If it is $20$ and $m=33,44,66$ the K3 surfaces are described in \cite{kondo1} and the fixed locus of the automorphism of order eleven is one point and one elliptic curve. If $m\not=33,44,66$ we argue in the same way as in the generic case. We now assume that the rank of the transcendental lattice is $10$. Since $\rank(NS(X))=12$ by \cite{niksympl} the only possibility is that $G_{0}$ is generated by a symplectic involution $\iota$. Observe that as before the case $m=22$ is not possible, hence $m=11$ and so $G$ has order $22$. If $G$ is a cyclic group of order $22$ then $\iota$ and $\eta$ commute and so the fixed locus of $\iota$ which are eight isolated fixed points must be permuted or fixed by $\eta$. This is not possible since $\eta$ has order  eleven and only two isolated fixed points. If $G$ is the dihedral group of order $22$ then the product $\eta\circ\iota$ has order two but the action on the holomorphic 2-form is the multiplication by an eleventh root of unity, which is not possible. Hence $G_0=\{1\}$ and $m=11$.    
\eprf
\begin{rem}{\rm 
In the paper \cite{OZ}, Oguiso and Zhang study the case of a K3 surface with a non--symplectic
automorphism of order eleven, they show the same result as Theorem \ref{711} for $p=11$. We gave here a unitary proof for both $p=7$ and $p=11$.}
\end{rem}

\section{K3 surfaces with a symplectic and non--symplectic automorphisms of the same order}\label{section: K3 surfaecs with sympl and non sympl}

Here we recall some basic facts on symplectic and non--symplectic automorphisms. One important point is that one can associate some lattices to an automorphism with a given fixed locus.

\begin{defi}{\rm (\cite[Definition 4.6]{niksympl})} We say that $G$ has a unique action on $\Lambda_{K3}$ if, given two embeddings $i: G\hookrightarrow Aut(X)$, $i': G\hookrightarrow Aut(X')$ such that $G$ is a group of symplectic automorphisms on the K3 surfaces $X$ and $X'$, there exists an isometry $\phi:H^2(X,\Z)\ra H^2(X',\Z)$ such that $i'(g)^*=\phi\circ i(g)\circ \phi^{-1}$ for all $g\in G$.\end{defi}

\begin{theorem}\label{theorem: unique action}{\rm (\cite[Theorem 4.7]{niksympl})} Let $G$ be a finite abelian group acting symplectically on a K3 surface. Then $G$ has a unique action on $\Lambda_{K3}$, hence the lattice $\Omega_G:=(\Lambda_{K3}^G)^{\perp}$ is uniquely determined by $G$, up to isometry.\end{theorem}

\begin{theorem}\label{theorem: X admits G iff omegaG in NS(X)}{\rm (\cite[Theorem 4.15]{niksympl}) } Let $G$ be a finite abelian group. A K3 surface $X$ admits $G$ as group of symplectic automorphisms if and only if the lattice $\Omega_G$ is primitively embedded in $NS(X)$.\end{theorem}
The lattices $\Omega_G$ are computed in \cite{aliale1}, \cite{aliale2} for each abelian group $G$. 
The fixed locus of a symplectic automorphism of prime order does not depend on the K3 surface on which the automorphism acts and it consists of a finite number of points. The situation is different for the non--symplectic automorphisms. The possible fixed loci of a non--symplectic automorphism of prime order on a K3 surface are listed in \cite{niknonsympl} in the case of involutions and in \cite{michialet} in the other cases.

\begin{defi}\label{elementary} Let $l$ be a prime number. A lattice is called $l$-elementary if its discriminant group is $(\Z/l\Z)^a$ for a certain non negative integer $a$. We then call $a:=l(L)$ the length of $L$.\\
Let $L$ be a 2-elementary lattice, then one defines the invariant $\delta$ of $L$ to be equal to $0$ if the discriminant form of $L$ takes values in $\Z$ and equal to $1$ otherwise.\end{defi}
We recall the following results on $l$-elementary lattices :
\begin{theorem}\label{niki}{\rm (\cite[Theorem 3.6.2]{nikbilinear})} An indefinite even $2$-elementary lattice is uniquely determined by $\delta$, its signature and the integer $a$.  
\end{theorem}
This can be generalized by the following theorem which gives also special results in the hyperbolic case:

\begin{theorem}\label{rudakov}{\rm (\cite{RS},\cite[Theorem 1.1]{michialet}) }
An even, indefinite, $p$-elementary lattice of rank $r$ for $p \neq 2$ and $r\geq 2$ is 
uniquely determined by its signature and the integer $a$.

For $p\neq 2$ a hyperbolic $p$-elementary lattice with invariants $a,r$  exists if and only if
the following conditions are satisfied:
$a\leq r$, $r\equiv 0 \pmod 2$ and 
\begin{equation*}
\begin{cases}
\text{for } a\equiv 0 \pmod 2, & r\equiv 2 \pmod 4 \\
\text{for } a\equiv 1 \pmod 2, & p\equiv (-1)^{r/2-1} \pmod 4
\end{cases}
\end{equation*}
and $r>a>0$, if $r\not \equiv 2 \pmod 8$.
\end{theorem}
Finally we recall some results on non--symplectic automorphisms:

\begin{theorem}{\rm (\cite{niksympl}, \cite{niknonsympl}, \cite{michiale}, \cite{michialet})}
Assume that a K3 surface $X$ has a non-sympletic automorphism $\eta$ of finite order $m$. Then we have :
\begin{itemize}
\item[1)] The Euler function $\varphi(m)$ divides $\rank(T_X)$. 
\item[2)] The lattice $H^2(X,\Z)^{\eta}$ is a hyperbolic lattice and is primitively embedded in $NS(X)$.
\item[3)] Let the order of $\eta$ be a prime number $p$, then the lattice $H^2(X,\Z)^{\eta}$ is a $p$-elementary lattice and the fixed locus $\Fix(\eta)$ determines uniquely the invariants $r$ and $a$ of the lattice $H^2(X,\Z)^{\eta}$.\end{itemize}
\end{theorem}
As a consequence of the previous results one obtains that if $X$ is a K3 surface admitting a non--symplectic automorphism $\eta$ of prime order $p$, the fixed locus determines $H^2(X,\Z)^{\eta}$ uniquely if $p$ is odd and gives at most two possible choices for $H^2(X,\Z)^{\eta}$ if $p=2$. Here we want to analyze K3 surfaces, $X$, admitting both a non--symplectic automorphism and a symplectic automorphism of the same order. This order is at most 8, in fact if $\Z/n\Z$ is a group of symplectic automorphisms on a complex K3 surface, then $n\leq 8$ (cf. \cite{niksympl}). Since $X$ admits $G:=\langle\sigma \rangle$ as group of symplectic automorphism, we have $\Omega_G\subset NS(X)$ (cf. Theorem \ref{theorem: X admits G iff omegaG in NS(X)}). The lattice $\Omega_G$, characterizing the K3 surfaces admitting $G$ as group of symplectic automorphisms, is negative definite. The K3 surfaces with non-symplectic automorphism of finite order are always algebraic (cf. \cite[Theorem 3.1]{niksympl}), hence the N\'eron--Severi group of $X$ contains a class with a positive self-intersection and so $\rank(NS(X))\geq \rank (\Omega_G)+1$. These observations lead us to Table \ref{tablerank}, where the possible rank of the N\'eron-Severi group and of the transcendental lattice for a K3 surface having a symplectic and a non-symplectic automorphism of the same order $m$ are given. The rank of the transcendental lattice determines the number of moduli of the family of K3 surfaces : 
\begin{eqnarray}\label{tablerank} 
\begin{array}{|c|c|c|c|} 
\hline
m&\rank (NS(X))&\rank (T_X)& \mbox{moduli}\\
\hline
2&\geq 9&\leq 13& \leq 11\\
\hline
3&14,~16,~18,~20& 8,~6,~4,~2& 3,~2,~1,~0\\
\hline
4&16,~18,~20&6,~4,~2&2,~1,~0\\
\hline
5&18&4&0\\
\hline
6&18,20&4,~2&1,~0 \\
\hline
7&-&-&-\\
\hline
8&-&-&-\\
\hline
\end{array}
\end{eqnarray}

We recall the contruction of a space parametrizing K3 surfaces with non-symplectic automorphism $\eta$ of prime order $p$ (cf. \cite{dolgakondo}, \cite{kondo}). 
We denote by $\eta^*$ the operation of $\eta$ on $H^2(X,\Z)$. By Nikulin (cf. \cite[Theorem 3.1]{niksympl}) the eigenvalues of $\eta$ on $T_X\otimes \mathbb{C}$  are $\zeta,\ldots,\zeta^{p-1}$, 
where $\zeta$ is a primitive $p$-root of unity. The value $1$ is not an eigenvalue of $\eta$. Hence one has a decomposition in eigenspaces:
$$
T_X\otimes \mathbb{C}=T_{\zeta}\oplus \ldots\oplus T_{\zeta^{p-1}}.
$$
Put
$$
\mathcal{B}=\{z\in\mathbb{P}(T_{\zeta})~:~(z,z)=0,~(z,\bar{z})>0\}
$$
and 
$$
\Gamma=\{\gamma\in O(T_X)~:~\gamma\circ\eta^*=\eta^*\circ\gamma\}.
$$
It is easy to see that $\mathcal{B}$ is a complex ball of dimesion $(\rank (T_X)/(p-1))-1$ for $p\geq 3$ and it is a type IV Hermitian symmetric space of dimension $\rank(T_X)-2$ if $p=2$. Then the generic point of $\mathcal{B}/\Gamma$ corresponds to a K3 surface with a non--symplectic automorphism $\eta'$ of order $p$ and $\eta'^*=\eta^*$. There is a birational map from $\mathcal{B}/\Gamma$ to the moduli space of K3 surfaces with a non--symplectic automorphism of order $p$ (see \cite{dolgakondo} and \cite{michialet} for a more detailed description). A similar construction holds also if the order of $\eta$ is  not a prime number. One must consider the decomposition of $T_X\otimes \mathbb{C}$ into eigenspaces corresponding to primitive $m$--roots of the unity and  the dimension of the moduli space is $(\rank (T_X)/\varphi(m))-1$. From these remarks it follows:

\begin{theorem}\label{theorem: no sympl e non for p=7,8}
There are no K3 surfaces having both a symplectic and a non--symplectic automorphism of the same order $n$ if $n=7,8$ and there are at most countable many K3 surfaces having a symplectic and a non--symplectic automorphism of order 5. 
\end{theorem}
\bprf
The first assertion follows from the Table \ref{tablerank} and the last one is a consequence of the structure of the moduli space, which in this case is 0-dimensional.
\eprf

\begin{prop}\label{prop:sympl and non commute}
Let $ \eta$ be a non--symplectic automorphism of finite order $m$ on a K3 surface $X$ such that the action of $ \eta$ on $NS(X)$ is the identity. Let $\sigma$ be a symplectic automorphism of finite order $n$ on the same K3 surface $X$. Then $ \eta$ and $\sigma$ commute and $ \eta\circ \sigma$ is an automorphism of order lcm$(m,n)$ acting on the period of the K3 surface as $\eta$.
\end{prop}
\bprf
By assumption $\eta$ acts as the identity on the N\'eron--Severi group $NS(X)$ and $\sigma$ acts as the identity on the transcendental lattice $T_X$ (cf. \cite{niksympl}). Hence the actions of $ \eta^*$ and $\sigma^*$ on $NS(X)\oplus T_X$ (and so on $H^2(X,\Z)$) commute. Since $ \eta^*\circ\sigma^*$ is an effective Hodge isometry on $H^2(X,\Z)$, by the global Torelli theorem, it is induced by a unique automorphism on $X$ (which is $ \eta\circ \sigma$). But $ \eta^*$ and $\sigma^*$ commute, hence $ \eta^*\circ\sigma^*$=$\sigma^*\circ \eta^*$ which is induced by a unique automorphism on $X$, i.e. by $\sigma\circ \eta$. So $ \eta\circ \sigma=\sigma\circ  \eta$.\\
Let $\omega_X$ be the period of the K3 surface $X$. So $( \eta\circ\sigma)^*(\omega_X)= \eta^*(\omega_X)=\zeta_m\omega_X$ where $\zeta_m$ is a primitive $m$-root of unity. Moreover since $ \eta$ and $\sigma$ commute, $ \eta\circ\sigma$ is an automorphism of order lcm$(m,n)$. 
\eprf

With respect to the fixed locus we have the following trivial property:
\begin{lemma}\label{lemma: fixed loci sympl and non commuting}
If two automorphisms $\eta$ and $\sigma$ commute, then $ \eta(\Fix(\sigma))\subset \Fix(\sigma)$ and viceversa $\sigma(\Fix( \eta))\subset \Fix( \eta)$. In particular if $ \eta$ and $\sigma$ are two automorphisms on a K3 surface as in Proposition \ref{prop:sympl and non commute}, then the fixed points of $\sigma$ are either fixed or permuted by $ \eta$.
\end{lemma}

\section{Order two}\label{section: order 2}

We recall that $\Omega_{\Z/2\Z}\simeq E_8(2)$. Moreover we recall that the invariant lattice $S_{\eta}$ of a non-symplectic involution $\eta$ is identified by three invariants $(r,a,\delta)$, where $r:=\rank(S_\eta)$, $a:=l(S_{\eta})$ and $\delta\in\{0,1\}$. We call $S_{(r,a,\delta)}$ the 2-elementary lattice with signature $(1,r-1)$ and invariants $(r,a,\delta)$. If a non-symplectic involution $\eta$ acts trivially on the N\'eron--Severi group of a K3 surface, then the N\'eron--Severi group coincides with some $S_{(r,a,\delta)}$ and hence $r$ is the Picard number of the surface. 
In the following $X_{(r,a,\delta)}$ will be a K3 surface with a non--symplectic involution $\eta$ acting trivially on the N\'eron--Severi group and such that $NS(X_{(r,a,\delta)})\simeq S_{(r,a,\delta)}$.

\begin{prop}\label{prop: Xrad sympl then a geq 16-r}
Assume that $X_{(r,a,\delta)}$ admits also a symplectic involution, then
\begin{itemize}
\item[i)] $r\geq 9$,
\item[ii)]$a\geq 16-r$,
\item[iii)] if $a=16-r$, then $r=10$, $a=6$, $\delta=0$. 
\end{itemize}\end{prop}

\bprf By \cite{niksympl}, if  $r=\rho(X_{(r,a,\delta)})\leq 8$, then $X_{(r,a,\delta)}$ does not admit a symplectic involution.\\ In order to prove ii) we can assume $r>8$ (by i)) and $r\leq 16$ (otherwise the implication is trivial, because $a$ is a non negative integer).
Since $X_{(r,a,\delta)}$ admits a symplectic automorphism, there exists a primitive embedding $\varphi:E_8(2)\ra NS(X_{(r,a,\delta)})$. Let $L:=\varphi(E_8(2))^{\perp_{NS(X_{(r,a,\delta)})}}$, so $NS(X_{(r,a,\delta)})$ is an overlattice of finite index of $E_8(2)\oplus L$ (we identify $E_8(2)$ with its immage $\varphi(E_8(2))$). The rank of $L$ is $r-8$. The length of $E_8(2)\oplus L$ is $8+l(L)$. To obtain an overlattice of $E_8(2)\oplus L$ such that $E_8(2)$ is primitively embedded in it, one has to add divisible classes of type $(e_i+l_i)/2$, where $e_i\in E_8(2)$, $l_i\in L$, $e_i/2\in E_8(2)^{\vee}/E_8(2)$ and $l_i/2\in L^{\vee}/L$. If one adds a divisible class to the lattice $E_8(2)\oplus L$ the length of the lattice decreases by two. The maximal number of divisible classes we can add is $l(L)$ (in fact since $r\leq 16$ we get $\rank (L)\leq 8$, so $l(L)\leq 8$), so the minimal possible length of an overlattice of $E_8(2)\oplus L$ such that $E_8(2)$ is primitively embedded in it, is $8+l(L)-2(l(L))=8-l(L)$. The length of $L$ cannot be greater then the rank of $L$, hence $8-l(L)\geq 8-\rank(L)=8-(r-8)=16-r$. Hence if $X_{(r,a,\delta)}$ admits a symplectic involution, $r\geq 9$ and $a\geq 16-r$.\\
To prove iii) we can assume $r\leq 16$ (otherwise $a$ is negative but there exists no K3 surface $X_{(r,a,\delta)}$ with a negative value of $a$). Let $a=16-r$. Then the lattice $L$ is such that rank$(L)=l(L)$, so we can choose a basis of $L$, $\{l_i\}_{i=1,\ldots r-8}$, such that $\{l_i/n_i\}$ are generators of the discriminant group $L^{\vee}/L$. Moreover we are adding exactly $r-8$ classes to $E_8(2)\oplus L$ to obtain $NS(X_{(r,a,\delta)})$ (because $a=16-r$, so it is the minimal possible). Hence there are $r-8$ elements in $E_8(2)$, called $e_i$, $i=1,\ldots r-8$, such that $$v_i:=(e_i+l_i)/2\in NS(X_{(r,a,\delta)}).$$ Since $v_i\in NS(X_{(r,a,\delta)})$, we have $v_i^2\in 2\Z$ and $v_iv_j\in \Z$. Recalling that $e_i^2\in 4\Z$ and $e_ie_j\in 2\Z$ we obtain
$$v_i^2=\frac{e_i^2+l_i^2}{4}\in 2\Z\mbox{ implies }l_i^2\in4\Z\mbox{ and }v_iv_j=\frac{e_ie_j+l_il_j}{4}\in \Z\mbox{ implies }l_il_j\in2\Z$$ and so $L=M(2)$ for a certain even lattice $M$. Since $NS(X_{(r,a,\delta)})$ is a 2-elementary lattice of length $a$, $|d(NS(X_{r,a,\delta}))|=2^a=2^{16-r}$. The discriminant of the overlattice of $E_8(2)\oplus L\simeq E_8(2)\oplus M(2)$ obtained adding the $r-8$ classes $(e_i+l_i)/2$ is $(2^8\cdot 2^{r-8}d(M))/2^{2(r-8)}=2^{16-r}d(M)$. So $|d(M)|=1$, i.e. $M$ is a unimodular lattice. Moreover it is clear that $M$ has signature $(1,r-9)$. This, together with the condition $r\leq 16$, implies that $M\simeq U$, $r=10$ and hence $a=6$. Let us now compute the discriminant form of $NS(X_{(r,a,\delta)})$. We recall that it is an overlattice of $(E_8\oplus U)(2)$. The discriminant group of $NS(X_{(r,a,\delta)})$ is generated by a linear combination, with integer coefficients, of the elements generating the discriminant group of $(E_8\oplus U)(2)$. Since the discriminant quadratic form of $(E_8\oplus U)(2)$ takes values in $\Z$, the discriminant quadratic form of $NS(X_{(r,a,\delta)})$ takes value in $\Z$ too, so $\delta=0$.\eprf

\begin{prop}\label{prop: if Xra symplectic then Xr+1a+1 and Xr+1a-1 symplectic}
i) If $X_{(r,a,1)}$ admits a symplectic involution and the surface $X_{(r+1,a+1,1)}$ exists, then it admits also a symplectic involution. \\
ii) If $X_{(r,a,1)}$ admits a symplectic involution and the surface $X_{(r+1,a-1,1)}$ exists, then it admits also a symplectic involution. 
\end{prop}
\bprf
A K3 surface $W$ admits a symplectic involution if and only if $E_8(2)$ is primitively embedded in $NS(W)$ or, equivalently, if and only if its transcendental lattice $T_W$ is primitively embedded in $E_8(2)\oplus U\oplus U\oplus U\simeq E_8(2)^{\perp_{\Lambda_{K3}}}$. By the assumptions on $X_{(r,a,1)}$, $E_8(2)$ is primitively embedded in $NS(X_{(r,a,1)})\simeq S_{(r,a,1)}$.\\
The lattice $S_{(r,a,1)}\oplus \langle -2\rangle$ is an even 2-elementary lattice with invariant $r'=r+1$, $a'=a+1$, $\delta=1$ and its signature is $(1,r'-1)$. These data identify uniquely its isometry class and so we have that $S_{(r,a,1)}\oplus\langle -2\rangle\simeq S_{(r+1,a+1,1)}$. Since $E_8(2)$ is primitively embedded in $S_{(r,a,1)}$, it is also primitively embedded in $S_{(r+1,a+1,1)}$. This proves that $X_{(r+1,a+1,1)}$ admits a symplectic automorphism.\\
Let $T_{(r,a,1)}:= T_{X_{(r,a,1)}}$ (observe that it is well defined since, by the Theorem \ref{niki}, $T_{(r,a,1)}$ is uniquely determined by $(r,a,1)$), we observe that $\rank(T_{(r,a,1)})=22-r$ and the length of $T_{(r,a,1)}$ is $a$. By the assumptions on $X_{(r,a,1)}$, the lattice $T_{(r,a,1)}$ is primitively embedded in $E_8(-2)\oplus U\oplus U\oplus U$. Let us consider the 2-elementary even lattice $T_{(r+1,a-1,1)}$, identified by the data $\rank (T_{(r+1,a-1,1)})=21-r$, $l(T_{(r+1,a-1,1)})=a-1$ , $\delta=1$ and signature $(2,19-r)$. It is clearly isometric to the transcendental lattice of the K3 surface $X_{(r+1,a-1,1)}$. We have $T_{(r,a,1)}\simeq T_{(r+1,a-1,1)}\oplus \langle -2\rangle$. Since $T_{(r,a,1)}$ is primitively embedded in $E_8(2)\oplus U\oplus U\oplus U$, the lattice $T_{(r+1,a-1,1)}$ is primitively embedded in $E_8(2)\oplus U\oplus U\oplus U$, hence $X_{(r+1,a-1,1)}$ admits a symplectic involution.\eprf

\begin{ex}\label{ex: order 2 sympl and non rho=9}{\bf The case $r=9$, $a=9$, $\delta=1$.}\\{\rm Let $X$ be a K3 surface with $NS(X)\simeq \Z L\oplus E_8(2)$,  and $L^2=2$. By \cite{bertale}, $X$ admits a symplectic involution $\iota$ (indeed there is a primitive embedding of $E_8(2)$ in $NS(X)$) and a projective model as double cover of $\mathbb{P}^2$ branched along the sextic $D_6:=V(f_6(x_1,x_2)+ax_0^2f_4(x_1,x_2)+bx_0^4f_2(x_1,x_2)+x_0^6)$, where the $f_d$ are homogeneous polynomials of degree $d$ in $x_1$, $x_2$. The symplectic involution is induced by the involution  $\iota_{\mathbb{P}^2}:(x_0:x_1:x_2)\mapsto (-x_0:x_1:x_2)$ of $\mathbb{P}^2$ (cf.\ \cite{bertale}). More precisely if the equation of $X$ is $u^2=f_6(x_1,x_2)+ax_0^2f_4(x_1,x_2)+bx_0^4f_2(x_1,x_2)+x_0^6=0$ in $W\mathbb{P}(3,1,1,1)$, then $\iota$ acts on $W\mathbb{P}(3,1,1,1)$ as $(u:x_0:x_1:x_2)\ra (-u:-x_0:x_1:x_2)$.
 The fixed locus of $\iota_{\mathbb{P}^2}$ is a line, which intersects the sextic in six points $p_i$, $i=1,\ldots, 6$, and of a point $q=(1:0:0)$. This induces eight fixed points on $X$.\\ 
The K3 surface $X$ admits also a non--symplectic automorphism $ \nu$ which is the covering involution (i.e. it is the restriction to $X$ of $(u:x_0:x_1:x_2)\ra (-u:x_0:x_1:x_2)$).\\ 
By construction $\iota$ acts as $-1$ on $E_8(2)\subset NS(X)$ and as $+1$ on the polarization $L$ and on the transcendental lattice. The covering involution $\nu$ preserves the class of the polarization and acts as $-1$ on its orthogonal complement. Let $\eta$ be the composition $\nu\circ \iota$, so it is induced on $X$ by $(u:x_0:x_1:x_2)\ra (u:-x_0:x_1:x_2)$. The fixed locus of $\eta$ is the genus 2 curve $u^2=f_6(x_1,x_2)$ (the curve in $X$ associated to $x_0=0$).\\ 
On the lattice $\Z L\oplus E_8(2)\oplus T_X\hookrightarrow H^2(X,\Z)$ we have the following action of these involutions:
$$\begin{array}{llll}
&\Z L\oplus &E_8(2)\oplus &T_X\\
\iota&+1&-1&+1\\
\nu &+1&-1&-1\\
\eta&+1&+1&-1
\end{array}$$
The non-symplectic automorphism $\eta$ acts trivially on the Picard group and hence the surface $X$ and the involution $\eta$ gives a model for the general member of the family of K3 surfaces with a non--symplectic involution $\eta$ with $H^2(X,\Z)^{\eta}\simeq \langle 2\rangle\oplus E_8(2)$ (or, which is the same, of the family of K3 surfaces with a non--symplectic involution with fixed locus a curve of genus 2).\\
We observe that $\Fix(\iota)\not \subset \Fix(\eta)$ and $\Fix(\eta)\not \subset \Fix(\iota)$.}\end{ex}

\begin{ex}\label{ex: order 2 sympl e non U2E82}{\bf The case $r=10$, $a=10$, $\delta=0$.} {\rm Let us consider a K3 surface $X$ such that $NS(X)\simeq U(2)\oplus E_8(2)$. It is well known that $X$ admits an Enriques involution, i.e. a fixed point free involution. It is also clear that $E_8(2)\hookrightarrow NS(X)$, hence this K3 surface admits also a symplectic involution.  In \cite[Chapter V, Section 23]{bpv} a 10-dimensional family of K3 surfaces admitting an Enriques involution is presented. Since the dimension of the family of K3 surfaces admitting an Enriques involution is 10 dimensional, the generic member of the family described in \cite[Chapter V, Section 23]{bpv} is the surface $X$.
It admits a 2:1 map to $\mathbb{P}^1\times \mathbb{P}^1$. Let us consider the involution $\iota_{\mathbb{P}^1\times\mathbb{P}^1}:((x_0:x_1);(y_0:y_1))\ra ((x_0:-x_1);(y_0:-y_1))$. It has four isolated fixed points $p_1=((0:1);(0:1))$, $p_2=((0:1);(1:0))$, $p_3=((1:0);(0:1))$, $p_4=((1:0);(1:0))$. Let $D_{4,4}$ be the curve \begin{eqnarray*}\begin{array}{c}ax_0^4y_0^4+bx_0^4y_0^2y_1^2+cx_0^4y_1^4+dx_0^2x_1^2y_0^4+ex_0^2x_1^2y_0^2y_1^2+fx_0^2x_1^2y_1^4+gx_1^4y_0^4+hx_1^4y_0^2y_1^2+lx_1^4y_1^4\\
+mx_0^3x_1y_0^3y_1+nx_0^3x_1y_0y_1^3+ox_0x_1^3y_0^3y_1+px_0x_1^3y_0y_1^3=0,
\end{array}
\end{eqnarray*} 
i.e. the invariant curve under the action of $\iota_{\mathbb{P}^1\times \mathbb{P}^1}$.  The double cover of $\mathbb{P}^1\times \mathbb{P}^1$ branched along $D_{4,4}$ is the K3 surface $X$. As in the previous example one obtains an equation of $X$ as $u^2=D_{4,4}$. The following three involutions act on $X$: 
$$\begin{array}{c} \iota:(u:(x_0,x_1):(y_0:y_1))\ra (u:\iota_{\mathbb{P}^1\times \mathbb{P}^1}((x_0:x_1):(y_0:y_1))),\\
\nu:(u:(x_0,x_1):(y_0:y_1))\ra (-u:((x_0:x_1):(y_0:y_1))),\\
\eta:(u:(x_0,x_1):(y_0:y_1))\ra (-u:\iota_{\mathbb{P}^1\times \mathbb{P}^1}((x_0:x_1):(y_0:y_1))),\end{array}$$
where $\iota$ is a symplectic automorphism fixing the eight points which are the inverse image (with respect to the double cover) of the four poits $p_i$, $\nu$ is the covering involution and $\eta$ is the Enriques involution. We observe that $\eta=\iota\circ \nu$.  
On the lattice $U(2)\oplus E_8(2)\oplus T_X\hookrightarrow H^2(X,\Z)$ we have the following action of these involutions:
$$\begin{array}{llll}
&U(2)\oplus &E_8(2)\oplus &T_X\\
\iota&+1&-1&+1\\
\nu &+1&-1&-1\\
\eta&+1&+1&-1
\end{array}$$
The non--symplectic involution $\eta$ fixes the lattice $U(2)\oplus E_8(2)\simeq NS(X)$.}
\end{ex}

\begin{ex}\label{ex: order 2 sympl e non UE82}{\bf The case $r=10$, $a=8$, $\delta=0$.} {\rm Let us consider the K3 surface $X$ admitting an elliptic fibration with Weierstrass equation $y^2=x^3+A(t^2)x+B(t^2)$. It is clear that it admits two non--symplectic involutions: $\nu: (x,y;t)\ra (x,-y;t)$ and $\eta:(x,y;t)\ra (x,y;-t)$. We observe that $\nu$ acts only on the fibers and $\eta$ only on the basis. The N\'eron--Severi group of $X$ is isometric to $U\oplus E_8(-2)$ and the Mordell--Weil lattice of the fibration is isometric to $E_8(2)$ (\cite{aliell}). The involution $\nu$ acts as $-1$ on each fiber and preserves the class of the fiber and the class of the zero section. In particular it acts as $+1$ on the copy of $U$ and as $-1$ on its orthogonal complement. The involution $\eta$ preserves the class of the fiber (it sends fibers to fibers) and the class of each section (because it acts only on the basis). So $\eta$ is a non--symplectic involution acting trivially on the N\'eron--Severi group. The composition $\eta\circ \nu$ is a symplectic automorphism acting as $-1$ on $E_8(2)$, i.e. on the Mordell--Weil lattice.}  \end{ex}

\begin{ex}\label{ex: order 2 sympl e non (10,6,0)}{\bf The case $r=10$, $a=6$, $\delta=0$.} {\rm  Let us consider the elliptic fibration with equation $y^2=x(x^2+a(t)x+b(t))$, with $deg(a(t))=4$ and $deg(b(t))=8$. It has 8 reducible fibers of type $I_2$ and its Mordell-Weil group is isometric to $\Z/2\Z$. Hence it admits a symplectic automorphism of order 2, which is the translation by the 2-torsion section. Its N\'eron--Severi is isometric to $U\oplus N$, where $N$ is the Nikulin lattice (cf. \cite[Proposition 4.2]{bertale}). We recall that the Nikulin lattice is a rank 8 even negative definite lattice and its discriminant form is the same as the discriminant form of $U(2)^3$. Hence the invariants of the N\'eron--Severi lattice are $r=10$, $a=6$, $\delta=0$, and it admits a symplectic involution. We observe that it clearly admits also a non--symplectic involution $\eta:(x,y,t)\ra (x,-y,t)$. It acts trivially on the N\'eron--Severi group and its fixed locus are two rational curves (the zero section and the 2-torsion section) and one curve of genus 3, the bisection $x^2+a(t)x+b(t)=0$, which is a $2:1$ cover of $\mathbb{P}^1$ branched in the eight zeros of $a(t)^2-4b(t)=0$. }
\end{ex}

\begin{ex}\label{ex: involution delta=0 2-torsion}{\bf Elliptic fibration with a 2-torsion section. }{\rm
In the following table we list certain elliptic K3 surfaces. Each of them is the generic member of a family of K3 surfaces with a non--symplectic involution associated to certain values of $(r,a,\delta)$ and the non--symplectic involution acting trivially on the N\'eron--Severi group is $\eta:(x,y,t)\ra (x,-y,t)$ (i.e. it acts as $-1$ on each smooth fiber of the fibration). Moreover each of these elliptic K3 surfaces admits a 2-torsion section (cf. \cite{shimada}), and hence a symplectic involution, which is the translation by this 2-torsion section (Example \ref{ex: order 2 sympl e non (10,6,0)} is a particular case of this construction). We will denote by $k$ the number of rational curves fixed by $\eta$ and by $g$ the genus of the non rational curve fixed by $\eta$ (the computation of $g$ and $k$ is similar to the one done in Example \ref{ex: order 2 sympl e non (10,6,0)}). 
\begin{eqnarray*}
\begin{array}{|c|c|c|c|c|}
\hline
(r,a,\delta)&\mbox{ singular fibers }&\mbox{ Mordell--Weil lattice}&k&g\\
\hline
(10,6,0)&8I_2+8I_1&\Z/2\Z&2&3\\
\hline
(14,4,0)&III^*+5I_2+6I_1&\Z/2\Z&5&2\\
\hline
(14,6,0)&I_2^*+6I_2+4I_1&\Z/2\Z&4&1\\
\hline
(18,0,0)&I_{12}^*+6I_1&\Z/2\Z&9&2\\
\hline
(18,2,0)&2III^*+2I_2+2I_1&\Z/2\Z&8&1\\
\hline
(18,4,0)&4I_0^*&\Z/2\Z&8&-\\
\hline
\end{array}
\end{eqnarray*}
In both the cases $(18,0,0)$ and $(18,2,0)$ the symplectic involution is in fact a Morrison--Nikulin involution, i.e. a sympletcic involution switching two copies of $E_8(1)$ in the N\'eron--Severi group (cf. \cite{morrison}, \cite{bertale}). The elliptic fibration given in case $(18,0,0)$ is described in details in \cite{clingherdoran}. The case $(18,2,0)$ is a particular member of the family described in Example \ref{ex: order 2 sympl e non U2E82}, some of its elliptic fibrations and involutions are described in \cite{matthiashulek}.
}\end{ex}

\begin{theorem}\label{theorem: syml and non iff, order 2} The K3 surface $X_{r,a,1}$ admits a symplectic involution if and only if $a>16-r$.\\
The K3 surface $X_{r,a,0}$ admits a symplectic involution if and only if either $a>16-r$ or $a=6$, $r=10$.\end{theorem}
\bprf By Proposition \ref{prop: Xrad sympl then a geq 16-r}, $a>16-r$ is a necessary condition to have a symplectic involution on $X_{r,a,1}$. By the Example \ref{ex: order 2 sympl and non rho=9} the surface $X_{9,9,1}$ admits a symplectic involution. By Proposition \ref{prop: if Xra symplectic then Xr+1a+1 and Xr+1a-1 symplectic} $ii)$, this implies that $X_{9+k,9-k,1}$, $k=0,\ldots 9$ admits a symplectic involution. By the Proposition \ref{prop: if Xra symplectic then Xr+1a+1 and Xr+1a-1 symplectic} $i)$, this implies that $X_{9+k,9+k,1}$, $k=0,1,2$ admits a symplectic involution and hence $X_{9+h,9-h,1}$, $X_{10+h,10-h,1}$, $X_{11+h,11-h,1}$,
$h=0,\ldots,9$ admit a symplectic involution.  This proves the first statement, because all the acceptable values of $r,a$ such that $a>16-r$ are of types $X_{9+k+h,9-k+h,1}$, $k=0,\ldots, 9$, $h=0,1,2$.\\
By Proposition \ref{prop: Xrad sympl then a geq 16-r} if $X_{r,a,0}$ admits a symplectic automorphism, then $a>16-r$ or $a=6$, $r=10$. The viceversa is proved in the examples \ref{ex: order 2 sympl e non U2E82}, \ref{ex: order 2 sympl e non UE82}, \ref{ex: order 2 sympl e non (10,6,0)}, \ref{ex: involution delta=0 2-torsion}.\eprf
 
\section{Order three}\label{section: order 3}
\begin{theorem}\label{theorem: order 3} Let $X$ be a K3 surface with a non--symplectic automorphism of order 3, $\eta$, acting trivially on the N\'eron--Severi group. If $rank(NS(X)):=\rho(X)<14$, then $X$ does not admit a symplectic automorphism of order 3.\\ The surface $X$ admits a symplectic automorphism of order 3 if and only if the fixed locus of $\eta$ consists $n$ points and $n-3$ curves and $n\geq 6$.\end{theorem}
\bprf The first statement follows from Table \ref{tablerank}. In \cite{michiale} it is proved that the families of K3 surfaces with a non--symplectic automorphism $\eta$ of order 3 are identified by the fixed locus of the automorphism, in particular by the pair $(n,k)$ where $n$ is the number of fixed points and $k$ is the number of fixed curves.  The families of K3 surfaces admitting a non--symplectic automorphism $\eta$ of order 3 such that the rank of the N\'eron--Severi group is greater than 14 are the families associated to the pairs $(n,n-3)$, $(m,m-2)$ with $6\leq n\leq 9$ and $6\leq m\leq 8$.\\
Moreover, it is proved that the generic member of a family of K3 surfaces admitting such an $\eta$ with fixed locus $(n,n-3)$, $6\leq n\leq 9$, has transcendental lattice $T_n$: 
$$T_6:=U\oplus U(3)\oplus A_2^2,\ \  T_7:=U\oplus U(3)\oplus A_2,\ \  T_8:=U\oplus U(3),\ \ T_9:=A_2(-1).$$ 
By Theorem \ref{theorem: X admits G iff omegaG in NS(X)} it sufficies to show that $\Omega_{\Z/3\Z}$ is primitively embedded in the N\'eron--Severi group of the K3 surfaces $X_n$ with transcendental lattice $T_n$. The lattice $\Omega_{\Z/3\Z}$ and its orthogonal complement in the K3 lattice $\Lambda_{K3}$ are computed in \cite{aliale1}. In particular it is proved that $\Omega_{\Z/3\Z}^{\perp}\simeq U\oplus U(3)^2\oplus A_2^2$ and $\Omega_{\Z/3\Z}=K_{12}(-2)$ (the Coxeter-Todd lattice of rank $12$).  Comparing the lattices $T_i$, $i=6,7,8,9$ and $\Omega_{\Z/3\Z}^{\perp}$, one notices that $$T_i\subset \Omega_{\Z/3\Z}^{\perp},\ \ \mbox{so} \ \ T_i^{\perp}\supset \Omega_{\Z/3\Z}.$$
Thus the N\'eron--Severi lattices of the surfaces admitting a non--symplectic automorphism with fixed locus $(6,3)$, $(7,4)$, $(8,5)$, $(9,6)$ contains the lattice $\Omega_{\Z/3\Z}$ and hence admits a symplectic automorphism of order 3.\\
The generic member of a family of K3 surfaces admitting such an $\eta$ with fixed locus $(m,m-2)$, $6\leq m\leq 8$, has transcendental lattice $T'_m:=U\oplus U\oplus A_2^{8-m}$. We observe that $T'_{m}$ is primitively embedded in $T'_{m-i}$, $i=1,2$.
We show that the lattice unimodular $T'_8\simeq U\oplus U$ is not primitively embedded in $U\oplus U(3)^2\oplus A_2^2\simeq (\Omega_{\Z/3\Z})^{\perp}$, indeed assume the contrary, then there exists a lattice $L$ such that $U\oplus U\oplus L\simeq (\Omega_{\Z/3\Z})^{\perp}$. This would implies that $L$ is a 3-elementary lattice of rank 6, signature $(1,5)$  and discriminant group equal to $(\Z/3\Z)^6$. By \ref{rudakov}, there exists no such a lattice, hence $T'_8$ is not primitively embedded in $\Omega_{\Z/3\Z}^{\perp}$. Since $T'_8$ is primitively embedded in $T'_{7}$ and $T'_{6}$, also $T'_7$ and $T'_6$ cannot be primitively embedded in $(\Omega_{\Z/3\Z})^{\perp}$. Thus $\Omega_{\Z/3\Z}$ is not primitively embedded in the N\'eron--Severi group of the generic K3 surface with a non--symplectic automorphism $\eta$ with fixed locus $(m,m-2)$, $6\leq m\leq 8$. \eprf 

In \cite{michiale} the generic member of the families of K3 surfaces admitting a non--symplectic automorphism $\eta$ of order 3 with fixed locus given by $n$ points and $n-3$ curves is described as an isotrivial elliptic fibration. In \cite{alibert} the generic member of such families with $6\leq n\leq 9$ is described by a different isotrivial elliptic fibration. Here we consider the description given in \cite{alibert} and we show that the generic member of this families also admits a symplectic automorphism of order 3. Moreover, since the non--symplectic automorphism and the symplectic automorphism commute we obtain that their composition is a non--symplectic automorphism of order three on the surface. We will analyze it later.

Let us consider now the K3 surface $S_f$ admitting an elliptic fibration with equation
\begin{align}\label{eq: Sf}S_f\,:\,y^2w=x^3+f_6^2(\tau)w^3, \,\,\ \ \tau\in\mathbb{C}\end{align}
where $f_6(\tau)$ is a polynomial of degree three with at most double zeros. In \cite{alibert} it is proved that these surfaces can be obtained by a quotient of the product surface $E_\zeta\times C_f$ by an automorphism of order three, where $C_f$ is a curve with equation $z^3=f_6(\tau)$ and $E_\zeta$ is the elliptic curve with equation $v^2=u^3+1$.\\
The surfaces $S_f$ clearly admits the non--symplectic automorphism of order three $$ \eta:(x:y:w;t)\mapsto (\zeta x:y:w;t).$$ 
By the equation of the elliptic fibration it is clear that it admits always the following sections $$s:\tau\mapsto (0:1:0;\tau),\  \ \ t_1:\tau\mapsto (0:f_6(\tau):1;\tau)\ \ \ t_2:\tau\mapsto (0:-f_6(\tau):1;\tau).$$
The sections $t_1$ and $t_2$ are 3-torsion sections, indeed for a fixed value $\overline\tau$ of $\tau$ they correspond to inflectional points of the elliptic curve $y^2w=x^3+f_6^2(\overline{\tau})w^3$. This implies that the K3 surface $S_f$ admits a symplectic automorphism $\sigma$ of order three induced by the translation by the section $t_1$ (cf. \cite{aliale2}). Hence these K3 surfaces admit both a non--symplectic automorphism and a symplectic automorphism of order three.

\subsection{Fixed locus $(6,3)$.}
Let us assume that $f_6(\tau)$ has six zeros of multiplicity one. Up to projectivity, one can assume that three of the zeros of $f_6(\tau)$ are in 0, 1 and $\infty$, so the elliptic fibration on $S_f$ is
$$\mathcal{E}_6:\ \ y^2=x^3+\tau^2(\tau-1)^2(\tau-\lambda_1)^2(\tau-\lambda_2)^2(\tau-\lambda_3)^2.$$
We observe that this family has three moduli. The singular fibers are six fibers of type $IV$ (i.e. three rational curves meeting in one point).  The sections $s$, $t_1$ and $t_2$ meet the singular fibers in different components.  Let us denote by $C_i^j$, $i=0,1,2$  the rational components of the $j$-th singular fibers with $C_0^j\cdot s=1$, $C^j_1\cdot t_1=1$, $C_2^j\cdot t_2=1$.\\
The non-symplectic automorphism $ \eta$ fixes the sections $s$, $t_1$, $t_2$ (since it fixes the base of the fibration). Moreover since the curves in the fixed locus of $ \eta$ are smooth and disjoint (cf.\ \cite{michiale}) $ \eta$ can not fix the components of the reducible fibers, because they meet the sections, but the curves $C_i^j$ are invariant under $ \eta$ and so $ \eta$ has two fixed points on each of them. One is the intersection between the curves $C_i^j$ and the sections, the other is the intersection point of the three rational curves $C_i^j$ in the same reducible fiber.\\
Summing up the fixed locus of $ \eta$ is made up of three rational curves ($s$, $t_1$, $t_2$) and six isolated points (the singular points of the six fibers of type $IV$). Hence the family of K3 surfaces admitting the fibration 
$\mathcal{E}_6$ is a (sub)family of the family of K3 surfaces admitting a non--symplectic automorphism of order three with fixed locus of type $(6,3)$. The dimension of that family is three (cf. \cite{michiale}) and the dimension of the family of K3 surfaces with equation $\mathcal{E}_6$ is 3 too, 
hence these two families coincide. In fact the moduli space $\mathcal M_{6,3}$ is irreducible as shown in \cite{michiale}.\\

The symplectic automorphism $\sigma$ acts as a translation by the three torsion section. In particular it preserves the fiber of the fibration and acts on the sections and on the components of the reducible fibers in the following way:
$$s\mapsto t_1\mapsto t_2,\ \ \ C_0^j\mapsto C_1^j\mapsto C_2^j,\ \ j=1,2,3,4,5,6.$$
Clearly, it fixes the singular point of the reducible fibers, and this is exactly $\Fix(\sigma)$ since a symplectic automorphism of order three on a K3 surface fixes exactly $6$ isolated points (cf. \cite{niksympl}).\\
In particular we observe that $\Fix(\sigma)\subset \Fix( \eta)$ and the isolated fixed points of $ \eta$ are exactly the same as the isolated fixed points of $\sigma$.\\

Let us now consider the automorphism $ \eta\circ\sigma$, which is a non--symplectic automorphism of order 3.\\
By  \cite[pp. 29--31]{silvtate}, $\sigma$ acts in the following way, if $x_P\neq 0$: $$\sigma(x_P:y_P:1;\tau)=\left(\frac{(y_P-f_6(\tau))^2-x_P^3}{x_P^2}:\left(\frac{y_P-f_6(\tau)}{x_P}\right)\left[x_P-\frac{(y_P-f_6(\tau))^2}{x_P^2}\right]-f_6(\tau):1;\tau\right).$$
If $x_P$ is zero, then we obtain the sections $s$, $t_1$  and $t_2$ and we described before the action on them. Finally we obtain $$ (\eta\circ\sigma)(x_P:y_P:1;\tau)=\left(\zeta\frac{(y_P-f_6(\tau))^2-x_P^3}{x_P^2}:\left(\frac{y_P-f_6(\tau)}{x_P}\right)\left[x_P-\frac{(y_P-f_6(\tau))^2}{x_P^2}\right]-f_6(\tau):1;\tau\right).$$
One can directly check that this automorphism fixes the curve $\mathcal{C}:\  \ x^3=-4f_6(\tau)^2$. This curve is a 3-section for the fibration. Since $\mathcal{C}$  is a 3:1 cover of $\mathbb{P}^1$ totally ramified  
over the zeros of $f_6(\tau)$, by the Riemann-Hurwitz formula we obtain that $g(\mathcal{C})=4$. No curves of the singular fibers are fixed by $ \eta\circ\sigma$ (indeed they are not invariant under $\sigma$ and they are invariant under $ \eta$),  so in the fixed locus of $ \eta\circ \sigma$ there is only one curve, i.e. $\mathcal{C}$, and this curve has to meet the singular fibers in their singular point ($\mathcal{C}$ can not meet one component of the fiber in a point which is not on the other components, because the components $C_i^j$ are not invariant under $ \eta\circ\sigma$). So the non--symplectic automorphism $ \eta\circ \sigma$ of order three has fixed locus of type $(0,1)$ and the fixed curve is of genus 4. In particular this implies that  the family of K3 surfaces admitting a non--symplectic automorphism of order 3 with fixed locus $(6,3)$ is a subfamily of the family of the K3 surfaces admitting a non-symplectic automorphism of order 3 with fixed locus $(0,1)$. This can be directly checked comparing the transcendental lattices of the generic member of these two families, c.f. \cite{michiale}. One can prove similarly that  $ \eta^2\circ\sigma$ is a non--symplectic automorphism of order 3 with fixed locus $(0,1)$ too.

\subsection{The other cases} In the following table we give the results obtained for different choices of the polynomial $f_6$ (they can be proved as in case $(6,3)$). In all the cases the curves fixed by $\eta$ are rational curves, the points fixed by $\sigma$ are isolated fixed points also of $\eta$ and the curve fixed by $\eta\circ\sigma$ is the trisection $x^3=-4f_6(\tau)^2$. 
\begin{align}\begin{array}{|c|c|c|c|}
\hline
f_{6}(\tau)&\mbox{ Singular fibers}&\mbox{Fixed locus }\eta&\mbox{Fixed locus }\eta\circ\sigma\\
\hline
\tau^2(\tau-1)^2(\tau-\lambda_1)^2(\tau-\lambda_2)^2(\tau-\lambda_3)&6IV&(6,3)&(0,1)\\
\hline
\tau^2(\tau-1)^2(\tau-\lambda_1)^2(\tau-\lambda_2)^2&4IV+IV^*&(7,4)&(1,1)\\
\hline
\tau^2(\tau-1)^4(\tau-\lambda_1)^4&2IV+2IV^*&(8,5)&(2,1)\\
\hline
\tau^4(\tau-1)^4&3IV^*&(9,6)&(3,1)\\
\hline
\end{array}\end{align}

\section{Order five}\label{section: order 5}
\begin{theorem}\label{theorem: order 5 non sympl and non}
The generic members of the families of K3 surfaces with a non--symplectic automorphism of order 5 cannot admit a symplectic automorphism of order 5. \end{theorem}
\bprf By Table \ref{tablerank} the the generic member of a family of K3 surfaces admitting both a symplectic and a non--symplectic  automorphism of order 5 has a the transcendental lattice of rank 4. By \cite{michialet} this implies that the transcendental lattice of the generic member of the family is isometric to $U\oplus H_5$. If there exists a primitive embedding of $U\oplus H_5$ in $U\oplus U(5)\oplus U(5)\simeq (\Omega_{\Z/5\Z})^{\perp}$, then there exists a rank 2 sublattice $M$ of $U\oplus U(5)\oplus U(5)$ such that $U\oplus H_5\oplus M$ is an overlattice of finite index of $U\oplus U(5)\oplus U(5)$. In particular the length $l$ of the discriminant group of $U\oplus H_5\oplus M$ has to be greater or equal to the length $m$ of $U\oplus U(5)\oplus U(5)$. The discriminant group of $U\oplus U(5)\oplus U(5)$ is $(\Z/5\Z)^4$, so $m=4$. The discriminant group of $U\oplus H_5$ is $\Z/5\Z$ and the length of $M$ is at most its rank, i.e. 2, so $l\leq 3$. Hence there exists no pimitive embedding of $U\oplus H_5$ in $U\oplus U(5)\oplus U(5)$ .\eprf
 
Clearly the previous proposition does not implies that there exist no K3 surfaces with both a non--symplectic and a symplectic automorphism of order 5, but that if there exists such a K3 surface, it is not the generic member of a family of K3 surfaces with a non--symplectic automorphism. The following example proves that there exist K3 surfaces with both symplectic and non--symplectic automorphism of order 5. By Table \ref{tablerank} we know that these K3 surfaces are rigid and the rank of their transcendental lattice is 4.  In the previous sections we constructed K3 surfaces admitting a non-symplectic automorphism of order $p=2,3$ and a symplectic automorphism of the same order commuting with the non-symplectic one. To construct the following example we again require that the non-symplectic automorphism and the symplectic automorphism commute.
 
 \begin{ex}\label{ex: rigid order 5} {\rm The general member of a family of K3 surfaces with a non--symplectic automorphism of order 5 with  four isolated fixed points as fixed locus, is a double cover of the plane ramified on a sextic $C$ in the family $\mathcal{C}$:
$$
a_1x_0^6+a_2x_0^3x_1x_2^2+a_3x_0^2x_1^3x_2+x_0(a_4x_1^5+a_5x_2^5)+a_6x_1^2x_2^4=0
$$
and the non-symplectic automorphism $\eta$ is induced by the automorphism of $\PP^2$
$$
\bar{\eta}(x_0:x_1:x_2)=(x_0:\zeta x_1:\zeta^2 x_2).
$$
This preserves each sextic in $\mathcal{C}$. The generic sextic is smooth and $\Fix(\bar\eta)=\{(1:0:0),(0:1:0),(0:0:1)\}$, since two points are on the sextic, $\Fix(\eta)$ are exactly four distinct points. The automorphisms of order 5 of $\mathbb{P}^2$ preserving the sextic and commuting with $\bar{\eta}$ are $(x_0:x_1:x_2)\mapsto (\zeta^a x_0:\zeta^b x_1:\zeta^c x_2)$.\\
To obtain a symplectic automorphism we consider those authomorphisms with $\zeta^a,\zeta^b,\zeta^c$ equal to permutations of $1,\zeta,\zeta^4$ leaving invariant some sextic in the family. The only possibility is $\bar\sigma:(x_0:x_1:x_2)\mapsto (x_0:\zeta x_1\zeta^4 x_2)$ and we obtain the family:
$$
x_0(a_1x_0^5+a_4x_1^5+a_5x_2^5)=0.
$$
The automorphisms of $\PP^2$ commuting with $\bar\sigma$ are only diagonal matrices, hence we obtain a number of parameters equal to zero, as expected.
In fact there are three possibilities for the branch sextic (not isomorphic to eachother):
$$
\begin{array}{l}
x_0(x_1^5+x_2^5+x_0^5)=0\\
x_0(x_1^5+x_2^5)=0\\
x_0(x_0^5+x_1^5)=0\\
x_0^6=0
\end{array}
$$
The last one is a multiple line, so the double cover is not a K3 surface.
The second case are six lines: five meeting at $(1:0:0)$ and one line not passing through this point. The third case are six lines meeting at $(0:0:1)$. The double cover in these cases is not a K3 surface (the singularity is not simple). There remains only one case. Hence the K3 surface $S$ is a double cover of $\mathbb{P}^2$ branched over the sextic 
$$ 
x_0(x_1^5+x_2^5+x_0^5)=0,
$$
which admits both a symplectic automorphism $\sigma$ induced by $\bar\sigma$ and a non--symplectic automorphism of order 5. The K3 surface $S$ admits also other non--symplectic automorphisms of order 5: \\
$\bullet$ $\nu$, induced by the automorphism  $\bar{\nu}:(x_0:x_1:x_2)\mapsto (x_0:\zeta x_1:x_2)$ of $\mathbb{P}^2$, it fixes a curve of genus 2 (the pullback of the line $x_1=0$) and a fixed point (the pullback of the point $(0:1:0)$);\\
$\bullet$ $\mu$ induced by $\bar{\mu}:(x_0:x_1:x_2)\mapsto (\zeta x_0:x_1:x_2)$, it fixes one rational curve (the pull back of the line $x_0=0$, which is contained in the branch locus) and 5 points (on the exceptional curve of the blow up of the singular points of the branch curves).
}\end{ex}
\begin{rem}{\rm The moduli space of K3 surfaces admitting a non--symplectic automorphism of order 5 has two irreducible components, corresponding to K3 surfaces with a non--symplectic automorphism of order 5 fixing only isolated points or fixing at least one curve. The K3 surface $S$ is in the intersection of these two components, indeed it admits a non--symplectic automorphism $\eta$ of order 5 fixing only isolated points and non--symplectic automorphisms $\nu$ and $\mu$,  of order 5, fixing at least one curve.}\end{rem}
\begin{rem}{\rm Observe that the symplectic automorphism $\sigma$ and the non--symplectic automorphism $\eta$ fix the same four points. This is a consequence of Lemma \ref{lemma: fixed loci sympl and non commuting}.}\end{rem}  

\section{Order four}\label{section: order 4}
By Table \ref{tablerank}, if a K3 surface admits both a symplectic and a non--symplectic automorphism of order 4, the rank of its transcendental lattice is $6,4,2$. Each of these transcendental lattices correspond to the generic member of a family of K3 surfaces admitting a non--symplectic automorphism of order 4, and the dimension of these families is respectively $2,1,0$. We denote these families by $\mathcal{M}_2$, $\mathcal{M}_1$ and $\mathcal{M}_0$ respectively. 
\begin{theorem}\label{theorem: no 2 dimensional family order 4} The generic member of the family $\mathcal{M}_2$ does not admit a symplectic automorphism of order 4.\end{theorem}
\bprf %By Table \ref{tablerank} the rank of the transcendental lattice of a K3 surface $S$ admitting both a non--symplectic and a symplectic automorphism of order 4 is at most 6. Hence the dimension of the moduli space is at most 2. 
Let $S$ be the generic member in the family $\mathcal{M}_2$ and $\eta$ be the non-symplectic automorhpism of order 4. It acts on the transcendental lattice with eigenvalues $i$ and $-i$ and on the N\'eron--Severi group with eigenvalues $1$ and $-1$  (cf. \cite{niksympl}). Thus $\eta^2$ acts on the transcendental lattice with eigenvalue $-1$ and on the N\'eron--Severi lattice with eigenvalue $1$. This implies that $T_S$ is a 2-elementary lattice. Suppose that $S$ admits also a symplectic automorphism of order $4$. Then $T_S\subset (\Omega_{\Z/4\Z})^{\perp}$. The discriminant group of $\Omega_{\Z/4\Z}^{\perp}$ is $(\Z/2\Z)^2\times (\Z/4\Z)^4$ and its rank is 8 (cf. \cite{aliale2}). If $T_S$ admits a primitive embedding in $\Omega_{\Z/4\Z}^{\perp}$, then there exists a rank 2 lattice $L$ such that $\Omega_{\Z/4\Z}^{\perp}$ is an overlattice of finite index of $T_S\oplus L$. The length of $\Omega_{\Z/4\Z}^\perp$ over $\Z/4\Z$ is 4, but the length of $T_S\oplus L$ over $\Z/4\Z$ is at most two, $T_S$ being a 2-elementary lattice and $L$ a lattice of rank 2. Hence it is impossible to find a lattice $L$ such that $\Omega_{\Z/4\Z}^{\perp}$ is an overlattice of finite index of $T_S\oplus L$, thus the general member of a 2-dimensional family of K3 surfaces admitting a non--symplectic automorphism of order 4 does not admit a symplectic automorphism of order 4.\eprf

\begin{prop}\label{prop: rigid auto 4} The K3 surface in the zero dimensional family of K3 surfaces admitting a non--symplectic automorphism of order 4, whose square fixes ten rational curves, admits a symplectic automorphism of order 4.\end{prop}
\bprf There exists only one family of K3 surface admitting a non--symplectic involution fixing 10 rational curves, this family is rigid and its transcendental lattice is $\langle 2\rangle^2$.   A model for this K3 surface is given in \cite{vinberg}: it is a $4:1$ cover of $\mathbb{P}^2$ branched along four lines, so up to a choice of coordinates on $\mathbb{P}^2$, it has an equation of type $$w^4=xy(x-z)(y-z).$$ The cover automorphism is the non--symplectic automorphism of order 4 whose square fixes ten rational curves (cf. \cite{aliCY}) . The automorphism of $\mathbb{P}^2$ given by $(x:y:z)\mapsto (x-z: x:x-y)$ preserves the branch locus of the $4:1$ cover (permuting the lines) and induces the automorphism of the surface $(w:x:y:z)\mapsto (w:x-z: x:x-y)$. It is a symplectic automorphism (this can be shown computing the action on the holomorphic 2--form).\eprf

\begin{rem}{\rm
The full automorphisms group of the surface in Proposition \ref{prop: rigid auto 4} is computed in \cite{vinberg}.} \end{rem}
\begin{ex}\label{ex: 1 dim order 4}{\rm Let us now consider the family of elliptic K3 surfaces: 
$$\mathcal{E}_4:\ \ wy^2=x^3+xw^2(at^8+bt^4+c).$$
We observe that this elliptic fibration has a 2-torsion section $t\mapsto (0:0:1;t)$ and 8 fibers of type $III$ (i.e. two tangent rational curves) over the zeros of $at^8+bt^4+c$. Moreover we observe that this elliptic fibration is isotrivial and the generic fiber is isomorphic to the elliptic curve $v^2=u^3+u$ which admits an automorphism of order 4 with two fixed points.\\
It is easy to see that one can assume that $a=c=1$ (it is enough to consider the transformation $(y,x)\mapsto (\lambda^3y,\lambda^2 x)$ and to divide by $\lambda^6$ to put one of the parameters equal to one and then to consider a projective transformation of $\mathbb{P}^1$ to put the other parameters equal to one), hence this family is 1--dimensional. 
It admits the dihedral group $\Dh_4$ of order 8 as group of symplectic automorphisms (cf. \cite{aliell}) which is generated by:
 $$\sigma_4:(x:y:w; t)\mapsto (-x:-iy: w; it)\ \ \ \varsigma_2:(x:y:w;t)\mapsto (\frac{x}{t^4}:-\frac{y}{t^6}:w;\frac{1}{t}).$$ 
In particular this implies that $NS(\mathcal{E}_4)$ contains primitively the lattice $U\oplus \Omega_{\Dh_4}$, where $U$ is generated by the classes of the fiber and of the zero section of the fibration, which are invariant classes for the action of $\sigma_4$ and $\varsigma_2$. Hence $\rank (NS(\mathcal{E}_4))\geq \rank (U)+\rank (\Omega_{\Z/4\Z})=2+15$ and $\rank (T_{X_4})\leq 5$. Since the family $X_4$ is a one dimensional family, we deduce that $\rank (T_{\mathcal{E}_4})=4$. 
The K3 surface $X_4$ clearly admits the non--symplectic automorphism of order four 
$$ \eta:(x:y:w;t)\mapsto (-x:iy:w;t).$$ 
The automorphisms $ \eta$ and $\sigma_4$ commute, hence their composition is a non--symplectic automorphism of order 4. The automorphism $\sigma_4$ fixes two points ($t=0$, $t=\infty$) on the base. The fibers over these points are smooth. It is easy to see that $\sigma_4$ has two fixed points on each of these two elliptic curves 
 and the automorphism $\sigma_4^2$ fixes four points on each of them. The automorphism $\varsigma_2$ has eight fixed points, four on the elliptic curve over $t=1$ and four on the elliptic curve over $t=-1$.\\
By using standard arguments on elliptic fibrations on K3 surfaces we can compute the fixed locus of the order 4 non--symplectic automorphisms obtained by the composition of the automorphisms described before. We resume the fixed loci in the following table. In particular $n_4$ is the number of the isolated fixed points of the automorphism of order 4, $k_4$ (resp. $k_2$) is the number of rational curves fixed by the automorphism of order 4 (resp. by its square), $g_4$ (resp. $g_2$) is the genus of the non rational curve(s) fixed by the automorphism of order 4 (resp. by its square), $r_2$ is the number of curves of genus $g_2$ fixed by the square of the automorphism. 

$$\begin{array}{|c|c|c|c|c|c|c|}
\hline
\mbox{automorphisms}&n_4&\mbox k_4&g_4&k_2&g_2&r_2\\
\hline
\eta&8&2&-&2&3&1\\
\hline
\eta\circ\sigma_4,\eta\circ\sigma_4^3&4&0&1&0&1&2\\
\hline
\eta\circ\sigma_4^2, \eta\circ\varsigma_2, \eta\circ\varsigma_2\circ \sigma_4^i,\ i=1,2&4&0&-&2&3&1\\
\hline
\end{array}$$
}\end{ex}

\section{Order six}\label{section: order 6}
We recall that a K3 surface admits a non--symplectic automorphism of order 6 if and only if it admits a non--symplectic automorphism of order 3 (cf. \cite{Dillies}). We have the following
\begin{theorem}\label{theorem: order 3 and 6 sympl} Let $X$ be a generic K3 surface with a non--symplectic automorphism of order 3. Then $X$ does not admit a symplectic automorphism of order 6.\end{theorem}
\bprf If $X$ admits a symplectic automorphism of order 6, $\sigma_6$,  then it admits also a symplectic automorphism,  $\sigma_6^2$, of order 3. By Theorem \ref{theorem: order 3}, $X$ must be in the families  of K3 surfaces admitting a non--symplectic automorphism of order 3 with fixed locus $(n,n-3)$, for $6\leq n\leq 9$. \\
Observe that to prove the theorem it is enough to show that $\Omega_{\Z/6\Z}$ is not contained in $U\oplus E_8^2\oplus A_2$, which is the N\'eron--Severi group of the generic K3 surface $S$ admitting a non--symplectic automorphism of order 3 with fixed locus $(9,6)$. Assume the contrary: if $S$ admits a symplectic automorphism of order 6, then $\Omega_{\Z/6\Z}$ is primitively embedded in $U\oplus E_8^2\oplus A_2$. Arguing (as in proof of Propostion \ref{prop: Xrad sympl then a geq 16-r}) on the length and on the rank of these two lattices one obtains  that there is no primitive embedding of $\Omega_{\Z/6\Z}$ in $U\oplus E_8^2\oplus A_2$. 
\eprf
It is however possible that there exist families of K3 surfaces admitting both a symplectic and a non--symplectic automorphism of order 6 and by Table \ref{tablerank} these families have at most dimension 1. Here we provide an example of such a family.
We recall that the moduli space of the K3 surfaces admitting a non--symplectic automorphism of order 3 (and hence of order 6) has 3 irreducible components $\mathcal{M}_0$, $\mathcal{M}_1$, $\mathcal{M}_2$, where $\mathcal{M}_i$ is the family of K3 surfaces admitting a non--symplectic automorphism of order 3 fixing $i$ curves (cf. \cite{michiale}). Moreover we recall that the rank of the transcendental lattice of a K3 surface admitting both a symplectic and a non--symplectic automorphism is at most 4. Here we construct a 1-dimensional family of K3 surfaces with both a symplectic and a non--symplectic automorphism of order 6. This family is a subfamily of $\mathcal{M}_1$ and $\mathcal{M}_2$.
\begin{ex}\label{ex: 1 dim order 6}{\rm 
Let $\mathcal{E}_6$ be the family of K3 surfaces
$$\mathcal{E}_6:\ \ y^2w=x^3+(at^{12}+bt^6+c)w^3.$$
We observe that this elliptic fibration has 12 fibers of type $II$ (i.e. a cuspidal rational curve) over the zeros of $(at^{12}+bt^6+c)$. Moreover we observe that this elliptic fibration is isotrivial and the generic fiber is isomorphic to the elliptic curve $v^2=u^3+1$ which admits an automorphism of order 3 with three fixed points. It is easy to see that one can assume that $a=c=1$. So this family depends on one parameter. In particular this implies that $\rank(T_{\mathcal{E}_6})=4$.\\
It admits the symplectic automorphisms (cf. \cite{aliell}) $$\sigma_6:(x:y:w; t)\mapsto (\zeta_6^2 x:-y: w; \zeta_6t),\ \ \ \varsigma_2:(x:y:w;t)\mapsto (\frac{x}{t^4}:-\frac{y}{t^6}:w;\frac{1}{t}).$$ 
Moreover it clearly admits the non--symplectic automorphisms $$ \eta:(x:y:w;t)\mapsto (\zeta_6^2 x,-y:w;t)$$ and $$\nu:(x:y:w;t)\ra (x:y:w;\zeta_6 t).$$ The automorphism $ \eta$ and $\sigma$ commute, hence their composition is a non--symplectic automorphisms of order 6, in particular $\nu=\eta^5\circ\sigma_6$. The automorphism $\sigma_6$, $\sigma_6^2$ and $\sigma_6^3$ fix each the points $t=0$, $t=\infty$ on the basis of the fibration. The fibers over these points are smooth and the automorphisms fix 2, 6 respectively 8 points on them. The automorphism $\varsigma_2$ has eight fixed points, four on the elliptic curve over $t=1$ and four on the elliptic curve over $t=-1$.\\

By using standard arguments on elliptic fibrations, we compute the fixed locus of the order 6 non--symplectic automorphisms obtained by the composition of the automorphisms described before. We resume the fixed loci in the table below. We denote by $n_6$ the number of the isolated fixed points of the automorphism of order 6, by $k_6$ (resp. $k_3$, $k_2$) the number of rational curves fixed by the automorphism of order 6 (resp. by its square, by its cube), by $g_6$ (resp. $g_3$, $g_2$) the genus of the non rational curve(s) fixed by the automorphism of order 6 (resp. by its square, by its cube), by $r_2$ the number of curves of genus $g_2$ fixed by the cube of the automorphism.
$$\begin{array}{|c|c|c|c|c|c|c|c|c|c|}
\hline
\mbox{automorphisms}&n_6&\mbox k_6&g_6&n_3&k_3&g_3&k_2&g_2&r_2\\
\hline
\eta,\ \ \nu\circ\sigma_6^5=\eta^5&12&1&-&0&1&5&1&10&1\\
\hline
\nu=\eta^5\circ\sigma_6,\ \ \eta\circ\sigma_6&3&0&1&3&0&1&0&1&2\\
\hline
\eta\circ\sigma_6^{2j},\ \ j=1,2,\ \ \nu\circ\sigma_6^i,\ \ i=1,3,4 &5&0&-&3&0&1&1&10&1\\
\hline
\eta\circ\sigma_6^3,\ \ \nu\circ\sigma_6^2,\ \ \varsigma_2\circ\eta,\ \ \varsigma_2\circ\eta^5,\ \ \varsigma_2\circ\nu,\ \ \varsigma_2\circ\nu^5&6&0&-&0&1&5&0&1&2\\
\hline
\end{array}$$
We observe that the automorphisms $\eta^j\circ \sigma^h$ for $h=1,\ldots, 5, \, j=2,3,4$ are not purely non--symplectic.
}
\end{ex}

\addcontentsline{toc}{section}{ \hspace{0.5ex} References}

\end{document}